\documentclass[11pt]{amsart}
 
\usepackage{amsmath,amsthm}
\usepackage{amssymb}
\usepackage{euscript}
\usepackage[frame,ps,matrix,arrow,curve,rotate]{xy} 

\numberwithin{equation}{subsection}

\newcommand{\sqsp}{\renewcommand{\baselinestretch}{1.2}\tiny\normalsize}

\newcommand{\relphantom[1]}{\mathrel{\phantom{#1}}}

\newtheorem{thm}[subsection]{Theorem}
\newtheorem{lemma}[subsection]{Lemma}
\newtheorem{prop}[subsection]{Proposition}
\newtheorem{cor}[subsection]{Corollary}

\theoremstyle{definition}
\newtheorem{definition}[subsection]{Definition}

\newcommand{\llbrack}{\lbrack \lbrack}  
\newcommand{\rrbrack}{\rbrack \rbrack}  
\newcommand{\Thetabar}{\bar{\Theta}}
\newcommand{\Thetatilde}{\tilde{\Theta}}

\newcommand{\thetatilde}{\tilde{\theta}}
\newcommand{\Psitilde}{\tilde{\Psi}}
\newcommand{\psitilde}{\tilde{\psi}}

\newcommand{\Ftilde}{\tilde{F}}

\newcommand{\fbar}{\bar{f}}
\newcommand{\ftilde}{\tilde{f}}

\newcommand{\sumprime}{\sideset{}{'}\sum}  
\newcommand{\lprod}{\dashv}  
\newcommand{\rprod}{\vdash}  

\DeclareMathOperator{\Id}{Id} 
  
\DeclareMathOperator{\Hom}{Hom}

\DeclareMathOperator{\Ob}{Ob}

\begin{document}
\title{Deformation theory of dialgebra morphisms}
\author{Donald Yau}

\begin{abstract}
An algebraic deformation theory of dialgebra morphisms is obtained.
\end{abstract}

\address{Department of Mathematics, The Ohio State University Newark, 1179 University Drive, Newark, OH 43055, USA}

\maketitle
\sqsp


\section{Introduction}
\label{sec:intro}

Dialgebras were introduced by Loday \cite{loday} in the study of periodicity phenomenon in algebraic $K$-theory.  As discussed in the introduction to the volume \cite{loday}, dialgebras arise naturally when one tries to construct a conjectural bicomplex, the analogue of the $(b,B)$-bicomplex for cyclic homology, that would give rise to algebraic $K$-theory.  In a dialgebra, there are two binary operations, $\lprod$ and $\rprod$, satisfying five associative style axioms.   Dialgebras are related to Leibniz algebras as associative algebras are related to Lie algebras.  Indeed, if one defines a bracket on a dialgebra by putting $\lbrack x, y \rbrack := x \lprod y - y \rprod x$, then one obtains a Leibniz algebra.  There is also an analogue of the universal enveloping algebra functor \cite{lp}.  Dialgebras, therefore, are of intrinsic algebraic interest.

The purpose of this paper is to study algebraic deformations of dialgebra morphisms, following the pattern established by Gerstenhaber \cite{ger1}.  The original deformation theory of associative algebras, as developed by Gerstenhaber \cite{ger1}, is closely related to Hochschild cohomology.  The relative version, the deformation theory of associative algebra morphisms, is studied by Gerstenhaber and Schack in a series of papers \cite{gs1,gs2,gs3}.  Deformations of Lie algebra morphisms have been studied by Nijenhuis and Richardson \cite{nr} and, more recently, by Fr\'eiger \cite{freiger}.  Since a Lie algebra is a Leibniz algebra in which the bracket is skew-symmetric, it should be possible to construct the deformation cohomology for a morphism of Leibniz algebras, following the approach in the Lie algebra case \cite{nr}.  Instead of Chevalley-Eilenberg cohomology, one would use Leibniz cohomology \cite{lp}.  More related to this paper, Majumdar and Mukherjee \cite{mm} worked out the absolute case, the deformation theory of dialgebras.

Although we deal with dialgebra morphisms in this paper, our approach does shed some new light into the classical case.  When studying deformations of an associative algebra morphism $\psi$, a major part is to show that the obstructions to extending a $2$-cocycle to a full-blown deformation are $3$-cocycles in the deformation complex $C^*(\psi,\psi)$ controlling the deformations of $\psi$.  To do this, the usual approach is to make use of a pre-Lie product \cite{ger0} on $C^*(\psi,\psi)$.  However, it is a highly non-trivial matter to establish directly the existence of such a structure on $C^*(\psi,\psi)$.  One can bypass this difficulty by invoking the powerful Cohomology Comparison Theorem (CCT) \cite{gs1,gs2,gs3}.  In particular, it allows one to ``pull" such a structure to $C^*(\psi,\psi)$ from the Hochschild cochain complex $C^*(\psi^!,\psi^!)$ of an auxiliary associative algebra $\psi^!$.

In showing that the corresponding obstructions are $3$-cocycles in the deformations of dialgebra morphisms, we take a direct approach.  In fact, once the appropriate obstructions $\Ob(\Theta_t)$ are identified, we compute $\delta \Ob(\Theta_t)$ explicitly and show that they are $0$.  We do not need a dialgebra version of the CCT, which is yet to be formulated and proved.  There are two advantages to our approach.  First, our direct calculation makes it clear as to why the obstructions are $3$-cocycles, without going through a dialgebra version of the CCT and an auxiliary dialgebra.  Second, simply by identifying the products $\lprod$ and $\rprod$, our argument that the obstructions are $3$-cocycles also applies to the classical case of an associative algebra morphism.  This gives an alternative route to the classical approach, free of the CCT.  Since the CCT is a very important result and has applications well beyond deformation theory, a direct argument that does not involve it makes the deformation theory of algebra morphisms simpler and more transparent.

The rest of this paper is organized as follows.  Section \ref{sec:prelim} is a preliminary section, in which we recall the basic definitions about dialgebras and their cohomology.  The deformation complex of a dialgebra morphism $\psi$ is constructed in Section \ref{sec:def complex}.  Deformations of $\psi$ and their infinitesimals are introduced in Section \ref{sec:deformation}.  It is observed that the infinitesimal is a $2$-cocycle in the deformation complex (Lemma \ref{lem:inf}).  Rigidity of deformations is studied in Section \ref{sec:eq rig}.  In particular, it is shown, as expected, that the vanishing of the second cohomology group $HY^2(\psi, \psi)$ implies that $\psi$ is rigid (Corollary \ref{cor:rigidity}).  The obstructions to extending a $2$-cocycle to a deformation are identified in Section \ref{sec:main}.  It is shown that such obstructions are $3$-cocycles (Lemma \ref{lem:main}) and that the simultaneous vanishing of their cohomology classes is equivalent to the existence of an extension to a deformation (Theorem \ref{thm:main}).  As an immediate consequence, the vanishing of $HY^3(\psi,\psi)$ implies that every $2$-cocycle can be realized as the infinitesimal of some deformation (Corollary \ref{cor:main}).


\section{Preliminaries on dialgebras}
\label{sec:prelim}

In this section, we briefly review some basic definitions and constructions about dialgebras.  The reader is referred to the references \cite{fra,loday} for details.

\subsection{Dialgebras}
\label{subsec:dialgebras}

Let $K$ be a field.  All tensor products and $\Hom$ are taken over $K$.  A \emph{dialgebra} $D$ over $K$ is a $K$-module equipped with two $K$-bilinear maps $\lprod$, $\rprod \colon D \otimes D \to D$, satisfying the following five axioms:
   \begin{equation}
   \label{eq:axioms}
   \begin{split}
   x \lprod (y \lprod z) &\buildrel 1 \over = (x \lprod y) \lprod z \buildrel 2 \over = x \lprod (y \rprod z), \\
   (x \rprod y) \lprod z &\buildrel 3 \over = x \rprod (y \lprod z), \\
   (x \lprod y) \rprod z &\buildrel 4 \over = x \rprod (y \rprod z) \buildrel 5 \over = (x \rprod y) \rprod z 
   \end{split}
   \end{equation}
for all elements $x, y, z \in D$.  The maps $\lprod$ and $\rprod$ are called the \emph{left product} and \emph{right product}, respectively.  A \emph{morphism} $\psi \colon D \to E$ of dialgebras is a $K$-linear map that respects the left and the right products.  A \emph{representation} of $D$ is a $K$-module $M$ equipped with two left actions, $\lprod, \rprod \colon D \otimes M \to M$, and two right actions, $\lprod, \rprod \colon M \otimes D \to M$, satisfying the fifteen axioms obtained from \eqref{eq:axioms} by choosing one of $x, y, z$ to be from $M$.  In particular, a dialgebra $D$ is a representation of itself, and if $\psi \colon D \to E$ is a dialgebra morphism, then $E$ is a representation of $D$ in an obvious way.

\subsection{Dialgebra cohomology}
\label{subsec:dialg coh}

For $n \geq 0$, let $Y_n$ denote the set of planar binary trees with $n + 1$ leaves, henceforth abbreviated to \emph{$n$-trees}.  The two $2$-trees are denoted by $\lbrack 1 2 \rbrack$ and $\lbrack 2 1 \rbrack$, and the five $3$-trees are denoted by $\lbrack 1 2 3 \rbrack$, $\lbrack 2 1 3 \rbrack$, $\lbrack 1 3 1 \rbrack$, $\lbrack 3 1 2 \rbrack$, and $\lbrack 3 2 1 \rbrack$ \cite[Appendix A]{loday}.  If $y$ is an $n$-tree and if $0 \leq i \leq n$, then $d_i y$ denotes the $(n-1)$-tree obtained from $y$ by deleting the $i$th leaf.  The maps $d_i$ satisfy the usual simplicial identities.

Let $D$ be a dialgebra over $K$ and $M$ be a representation of $D$.  Define the \emph{module of $n$-cochains} $CY^n(D,M)$ to be the $K$-module $\Hom(K\lbrack Y_n \rbrack \otimes D^{\otimes n}, M)$.  The coboundary $\delta \colon CY^n(D,M) \to CY^{n+1}(D,M)$ is defined to be $\delta = \sum_{i=0}^{n+1} (-1)^i \delta^i$, where
   \begin{multline*}
   (\delta^i f)(y \otimes (a_1, \ldots, a_{n+1})) \\
   = 
   \begin{cases} 
   a_1 \circ^y_0 f(d_0y \otimes (a_2, \ldots a_{n+1})) & \text{ if } i = 0, \\
   f(d_iy \otimes (\ldots, a_{i-1}, a_i \circ^y_i a_{i+1}, a_{i+2}, \ldots)) & \text{ if } 1 \leq i \leq n, \\
   f(d_{n+1}y \otimes (a_1, \ldots, a_n)) \circ^y_{n+1} a_{n+1} & \text{ if } i = n + 1.
   \end{cases}
   \end{multline*}
Here $\circ_i \colon Y_{n+1} \to \lbrace \lprod, \rprod \rbrace$ is a certain function \cite[p.\ 75]{fra}.  As $\delta \circ \delta = 0$, one defines the \emph{dialgebra cohomology $HY^*(D,M)$ of $D$ with coefficients in $M$} to be the $n$th cohomology group of $(CY^*(D,M), \delta)$.


\section{Deformation complex of a dialgebra morphism}
\label{sec:def complex}

Throughout the rest of this paper, $K$ denotes an arbitrary but fixed field, and  $\psi \colon D \to E$ denotes a morphism of dialgebras over $K$.  Regard $E$ as a representation of $D$ via $\psi$ wherever appropriate.

\subsection{Deformation complex}
\label{subsec:def complex}

Define the \emph{module of $n$-cochains of $\psi$} by
   \begin{equation}
   \label{eq:n-cochain psi}
   CY^n(\psi, \psi) \, = \, CY^n(D,D) \times CY^n(E,E) \times CY^{n-1}(D,E).
   \end{equation}   
The coboundary $\delta \colon CY^n(\psi,\psi) \to CY^{n+1}(\psi,\psi)$ is defined by the formula
   \begin{equation}
   \label{eq:coboundary psi}
   \delta(\xi;\, \pi;\, \varphi) \, = \, (\delta \xi;\, \delta \pi;\, \psi \xi - \pi \psi - \delta \varphi) 
   \end{equation}
for $(\xi; \pi; \varphi) \in CY^n(\psi,\psi)$.  Here $\psi \xi$ and $\pi \psi$ in $CY^n(D,E)$ are the push-forwards:
   \begin{align*}
   (\psi \xi)(y \otimes (a_1, \ldots, a_n)) & \,=\, \psi(\xi(y \otimes (a_1, \ldots, a_n))), \\
   (\pi \psi)(y \otimes (a_1, \ldots, a_n)) & \,=\, \pi(y \otimes (\psi(a_1), \ldots, \psi(a_n)))
   \end{align*}
for $y \in Y_n$ and $a_1, \ldots, a_n \in D$.  Note the similarity here with the associative case \cite[p.\ 155]{gs3}.

\begin{prop}
\label{prop:def complex}
$(CY^*(\psi,\psi), \delta)$ is a cochain complex, i.e.\ $\delta \circ \delta = 0$.
\end{prop}

\begin{proof}
The right-most component of $(\delta \circ \delta)(\xi; \pi; \varphi)$ is $(\psi(\delta \xi) - (\delta \pi)\psi - \delta(\psi \xi - \pi \psi - \delta \varphi))$.  To finish the proof, one checks by direct inspection that $\psi(\delta \xi) = \delta (\psi \xi)$ and $(\delta \pi) \psi = \delta(\pi \psi)$.
\end{proof}

The cochain complex $(CY^*(\psi,\psi), \delta)$ is called the \emph{deformation complex of $\psi$}.  Define the \emph{$n$th dialgebra cohomology of $\psi$} by
   \[
   HY^n(\psi, \psi) \,\buildrel \text{def} \over =\, H^n(CY^*(\psi,\psi), \delta).
   \]
This is related to the dialgebra cohomologies of $D$, $E$ (with self coefficients), and $D$ with coefficients in $E$ in the following way.

\begin{prop}
\label{prop2:def complex}
If $HY^n(D,D)$, $HY^n(E,E)$, and $HY^{n-1}(D,E)$ are all trivial, then so is $HY^n(\psi, \psi)$.
\end{prop}

\begin{proof}
Let $\alpha = (\xi; \pi; \varphi) \in CY^n(\psi,\psi)$ be an $n$-cocycle.  Then by definition \eqref{eq:coboundary psi} and the hypothesis, one has that $\xi = \delta \xi^\prime$ and $\pi = \delta \pi^\prime$ for some $(n-1)$-cochains $\xi^\prime \in CY^{n-1}(D,D)$, $\pi^\prime \in CY^{n-1}(E,E)$.  Since $\delta \alpha = 0$, we have that
   \[
   \begin{split}
   0 
   &\,=\, \psi \xi - \pi \psi - \delta \varphi \\
   &\,=\, \psi(\delta \xi^\prime) - (\delta \pi^\prime)\psi - \delta \varphi \\
   &\,=\, \delta(\psi \xi^\prime) - \delta(\pi^\prime \psi) - \delta \varphi \\
   &\,=\, \delta(\psi\xi^\prime - \pi^\prime \psi - \varphi),  
   \end{split}
   \]
i.e.\ $(\psi\xi^\prime - \pi^\prime \psi - \varphi)$ is an $(n-1)$-cocycle.  It follows from the hypothesis that $(\psi\xi^\prime - \pi^\prime \psi - \varphi) = \delta \varphi^\prime$ for some $(n-2)$-cochain $\varphi^\prime \in CY^{n-2}(D,E)$ and, hence, $\alpha = \delta(\xi^\prime; \pi^\prime; \varphi^\prime)$.
\end{proof}


\section{Deformation and infinitesimal}
\label{sec:deformation}

\subsection{Deformation}

Recall from \cite{mm} that a \emph{deformation} of a dialgebra $D$ over $K$ is a pair of power series, $f_t = \left(f^l_t = \sum_{n=0}^\infty F^l_n t^n;\, f^r_t = \sum_{n=0}^\infty F^r_n t^n \right)
   = \sum_{n=0}^\infty (F^l_n;\, F^r_n)t^n$, satisfying the following three conditions: (i) Each $F^*_n \colon D \otimes D \to D$ is a $K$-bilinear map.  (ii) $F^l_0 = \lprod$, $F^r_0 = \rprod$.  (iii) By extension to power series, the $K$-module $D \llbrack t \rrbrack$ equipped with the products $f^l_t$ and $f^r_t$ becomes a dialgebra, denoted $D_t$ or $D_t(f)$.  Sometimes $D_t$ itself is referred to as the deformation.  The pair $(F^l_n;\, F^r_n)$ determines and is determined by the $2$-cochain $F_n \in CY^2(D,D)$, where $F_n(y \otimes (a_1, a_2)) = F^*_n(a_1, a_2)$ with $* = l$ (resp.\ $* = r$) if $y = \lbrack 2 1 \rbrack$ (resp. $y = \lbrack 1 2 \rbrack$).  For this reason, we will often write $f_t$ as $\sum_{n=0}^\infty F_n t^n$.

\begin{definition}
\label{def:deformation psi}
Let $\psi \colon D \to E$ be a dialgebra morphism over $K$.  Define a \emph{deformation of $\psi$} to be a triple $\Theta_t = \left(f_{D,t};\, f_{E,t};\, \Psi_t\right)$ in which:
   \begin{itemize}
   \item $f_{D,t} = \sum_{n=0}^\infty F_{D,n}t^n$ is a deformation of $D$;
   \item $f_{E,t} = \sum_{n=0}^\infty F_{E,n}t^n$ is a deformation of $E$;
   \item $\Psi_t \colon D_t \to E_t$ is a dialgebra morphism of the form $\Psi_t = \sum_{n=0}^\infty \psi_n t^n$, where each $\psi_n \colon D \to E$ is a $K$-linear map and $\psi_0 = \psi$.
   \end{itemize}
\end{definition}

Since there is only one $1$-tree, each $\psi_n$ can be identified with a $1$-cochain in $CY^1(D,E)$.  In particular, each $\theta_n = (F_{D,n}; F_{E,n}; \psi_n)$ is a $2$-cochain in $CY^2(\psi,\psi)$.  We will often write a deformation $\Theta_t$ as a power series itself, $\Theta_t = \sum_{n=0}^\infty \theta_n t^n$.

\subsection{Infinitesimal}

\begin{definition}
\label{def:infinitesimal}
The linear coefficient, $\theta_1 = (F_{D,1}; F_{E,1}; \psi_1)$, is called the \emph{infinitesimal} of the deformation $\Theta_t$.
\end{definition}

Note that $F_{D,1}$ and $F_{E,1}$ are the infinitesimals of $D$ and $E$, respectively \cite[Definition 3.2]{mm}.

\begin{lemma}
\label{lem:inf}
The infinitesimal $\theta_1$ of a deformation $\Theta_t$ of $\psi$ is a $2$-cocycle in $CY^2(\psi,\psi)$.  More generally, if $\theta_i = 0$ for $i = 1, 2, \ldots, n$, then $\theta_{n+1}$ is a $2$-cocycle.
\end{lemma}

\begin{proof}
It is proved in \cite[Lemma 3.3]{mm} that $F_{D,1}$ is a $2$-cocycle in $CY^2(D,D)$.  The same remark applies to $F_{E,1}$.  To finish the proof, notice that the condition that $\Psi_t$ be a dialgebra morphism is equivalent to the equality
   \begin{equation}
   \label{eq:Psi}
   \Psi_t\left(f^*_{D,t}(a,b)\right) = f^*_{E,t}\left(\Psi_t(a), \Psi_t(b)\right)
   \end{equation}
for $* \in  \lbrace l, r \rbrace$ and $a, b \in D$.  This in turn is equivalent to the condition
   \begin{equation}
   \label{eq:Psi'}
   \sum_{n=0}^N \psi_n F^*_{D,N-n}(a, b) = \sum_{i + j + k = N} F^*_{E,i}(\psi_j(a), \psi_k(b))
   \end{equation}
for $* \in \lbrace l, r \rbrace$, $a, b \in D$, and $N \geq 0$.  When $N = 0$, this just says that $\psi$ preserves the left and the right products.  But when $N = 1$, \eqref{eq:Psi'} can be rewritten as
   \[
   \begin{split}
   0 \,=\, \psi F^l_{D,1}(a,b) &- F^l_{E,1}(\psi(a), \psi(b)) \\ 
   &- (\psi(a) \lprod \psi_1(b) - \psi_1(a \lprod b) + \psi_1(a) \lprod \psi(b)), \\
   0 \,=\, \psi F^r_{D,1}(a,b) &- F^r_{E,1}(\psi(a), \psi(b)) \\
   &- (\psi(a) \rprod \psi_1(b) - \psi_1(a \rprod b) + \psi_1(a) \rprod \psi(b)).
   \end{split}
   \]
This is equivalent to saying that the $2$-cochain $(\psi F_{D,1} - F_{E,1}\psi - \delta\psi_1) \in CY^2(D,E)$ is equal to $0$.  Therefore, we have $\delta \theta_1 = 0$, as desired.

The second assertion is proved similarly.
\end{proof}


\section{Equivalence and rigidity}
\label{sec:eq rig}

\subsection{Equivalence}

Let $D_t(f)$ and $D_t(\ftilde)$ be two deformations of $D$.  Recall from \cite[Definition 4.1]{mm} that a \emph{formal isomorphism} $\Phi_t \colon D_t(f) \to D_t(\ftilde)$ is a power series $\Phi_t = \sum_{n=0}^\infty \phi_n t^n$ in which each $\phi_n \colon D \to D$ is a $K$-linear map and $\phi_0 = \Id_D$ such that 
   \begin{equation}
   \label{eq:formal iso D}
   \ftilde^*_t(a, b) \,=\, \Phi_t f^*_t(\Phi_t^{-1}(a), \Phi_t^{-1}(b))
   \end{equation}
for all $a, b \in D$ and $* \in \lbrace l, r \rbrace$.  Two deformations $D_t(f)$ and $D_t(\ftilde)$ are \emph{equivalent} if and only if there exists a formal isomorphism $D_t(f) \to D_t(\ftilde)$.

\begin{definition}
\label{def:formal iso}
Let $\Theta_t = (f_{D,t}; f_{E,t}; \Psi_t)$ and $\Thetatilde_t = (\ftilde_{D,t}; \ftilde_{E,t}; \Psitilde_t)$ be two deformations of a dialgebra morphism $\psi \colon D \to E$.  A \emph{formal isomorphism} $\Phi_t \colon \Theta_t \to \Thetatilde_t$ is a pair $\Phi_t = (\Phi_{D,t}; \Phi_{E,t})$, where $\Phi_{D,t} \colon D_t(f_D) \to D_t(\ftilde_D)$ and $\Phi_{E,t} \colon E_t(f_E) \to E_t(\ftilde_E)$ are formal isomorphisms, such that    
   \begin{equation}
   \label{eq:formal iso}
   \Psitilde_t  \,=\, \Phi_{E,t} \Psi_t \Phi_{D,t}^{-1}.
   \end{equation}
Two deformations $\Theta_t$ and $\Thetatilde_t$ are \emph{equivalent} if and only if there exists a formal isomorphism $\Theta_t \to \Thetatilde_t$.
\end{definition}

Here is a simple but very useful observation.  Given only a deformation $\Theta_t$ and a pair of power series $\Phi_t = (\Phi_{D,t} = \sum \phi_{D,n}t^n; \Phi_{E,t} = \sum \phi_{E,n}t^n)$ as above, one can define a deformation $\Thetatilde_t$ using \eqref{eq:formal iso D} (for both $D$ and $E$) and \eqref{eq:formal iso}.  The resulting deformation $\Thetatilde_t$ is automatically equivalent to $\Theta_t$.

\begin{thm}
\label{thm:infinitesimal}
The infinitesimal of a deformation $\Theta_t$ of $\psi$ is a $2$-cocycle in $CY^2(\psi,\psi)$ whose cohomology class is determined by the equivalence class of $\Theta_t$. 
\end{thm}

\begin{proof}
In view of Lemma \ref{lem:inf}, it remains to show that if $\Phi_t \colon \Theta_t \to \Thetatilde_t$ is a formal isomorphism, then the $2$-cocycles $\theta_1$ and $\thetatilde_1$ differ by a $2$-coboundary.  Write $\Phi_t = (\Phi_{D,t} = \sum_{n=0}^\infty \phi_{D,n} t^n;\; \Phi_{E,t} = \sum_{n=0}^\infty \phi_{E,n} t^n)$ and $\Thetatilde_t = (\ftilde_{D,t} = \sum (\Ftilde^l_{D,n};\, \Ftilde^r_{D,n})t^n;\, \ftilde_{E,t} = \sum (\Ftilde^l_{E,n};\, \Ftilde^r_{E,n})t^n;\, \Psitilde_t = \sum \psitilde_n t^n)$.  It is shown in \cite[Proposition 4.3]{mm} that $\delta \phi_{*,1} = F_{*,1} - \Ftilde_{*,1}$ in $CY^2(*,*)$ where $* = D, E$.

To finish the proof, observe that the linear coefficients on both sides of \eqref{eq:formal iso} yield the equality
   \[
   \psi_1 - \psitilde_1 \,=\, \psi \phi_{D,1} - \phi_{E,1}\psi.
   \]
It follows that the $1$-cochain $\alpha = (\phi_{D,1}; \phi_{E,1}; 0) \in CY^1(\psi,\psi)$ satisfies 
$\delta \alpha = \theta_1 - \thetatilde_1$, as desired.
\end{proof}

\subsection{Rigidity}

\begin{definition}
\label{def:rigidity}
A dialgebra morphism $\psi$ is said to be \emph{rigid} if and only if every deformation of $\psi$ is equivalent to the trivial deformation $(F_{D,0}; F_{E,0}; \psi)$.
\end{definition}

\begin{thm}
\label{thm:rigidity}
Let $\Theta_t = \left(f_{D,t};\, f_{E,t};\, \Psi_t\right) = \sum_{n=0}^\infty \theta_n t^n$ be a deformation of $\psi$ in which $\theta_i = 0$ for $i = 1, \ldots, m$ and $\theta_{m+1}$ is a $2$-coboundary in $CY^2(\psi,\psi)$.  Then there exists a deformation $\Thetatilde_t = \sum_{n=0}^\infty \thetatilde_n t^n$ of $\psi$ and a formal isomorphism $\Phi_t \colon \Theta_t \to \Thetatilde_t$ such that:
   \begin{enumerate}
   \item $\Phi_t = (\Phi_{D,t} = 1_D + \xi t^{m+1};\, \Phi_{E,t} = 1_E + \pi t^{m+1})$ for some $\xi \in CY^1(D,D)$ and  $\pi \in CY^1(E,E)$;
   \item $\thetatilde_i = 0$ for $i = 1, \ldots, m+1$. 
   \end{enumerate}
\end{thm}

To prove this Theorem, we need the following observation.

\begin{lemma}
\label{lem:rigidity}
Let $\theta$ be a $2$-coboundary in $CY^2(\psi,\psi)$.  Then there exists a $1$-cochain $\beta \in CY^1(\psi,\psi)$ of the form $\beta = (\xi; \pi; 0)$ such that $\theta = \delta \beta$. 
\end{lemma}

\begin{proof}
The proof here is identical with that of the usual case of an associative algebra morphism \cite[page 156]{gs3}.  Indeed, direct inspection shows that $\delta(\xi; \pi; \varphi) = \delta(\xi; \pi + \delta \varphi; 0)$ for any $1$-cochain $(\xi; \pi; \varphi) \in CY^1(\psi,\psi)$.
\end{proof}

\begin{proof}[Proof of Theorem \ref{thm:rigidity}]
Using Lemma \ref{lem:rigidity}, write $\theta_{m+1}$ as $\delta\beta = (\delta \xi; \delta \pi; \psi \xi - \pi \psi)$ for some $1$-cochains $\xi \in CY^1(D,D)$, $\pi \in CY^1(E,E)$.  Define a pair of power series $\Phi_t = (\Phi_{D,t} = 1_D + \xi t^{m+1};\, \Phi_{E,t} = 1_E + \pi t^{m+1})$.  Then define a deformation $\Thetatilde_t = (\ftilde_{D,t}; \ftilde_{E,t}; \Psitilde_t) = \sum \thetatilde_n t^n$ using \eqref{eq:formal iso D} (for both $D$ and $E$) and \eqref{eq:formal iso}.  It is then automatic that $\Phi_t \colon \Theta_t \to \Thetatilde_t$ is a formal isomorphism and that condition (1) in Theorem \ref{thm:rigidity} holds.

To check condition (2), we compute modulo $t^{m+2}$:
   \begin{equation}
   \label{eq:rigidity psitilde}
   \begin{split}
   \Psitilde_t 
   &\,=\, \Phi_{E,t} \Psi_t \Phi_{D,t}^{-1} \\
   &\,\equiv\, (1_E + \pi t^{m+1})(\psi + \psi_{m+1}t^{m+1})(1_D - \xi t^{m+1}) \\
   &\,\equiv\, \psi + (\psi_{m+1} - \psi\xi + \pi\psi)t^{m+1} \\
   &\,=\, \psi.
   \end{split}
   \end{equation}
If we write $\ftilde_{*,t} = \sum_{n=0}^\infty \Ftilde_{*,n}t^n$, where $* = D, E$, then the proof that $\Ftilde_{*,i} = 0$ for $i = 1, \ldots, m+1$ is given in \cite[Theorem 4.5]{mm}.  Combined with \eqref{eq:rigidity psitilde}, we conclude that condition (2) holds as well.
\end{proof}

Applying Lemma \ref{lem:inf} and Theorem \ref{thm:rigidity} repeatedly, we obtain the following cohomological criterion for rigidity.

\begin{cor}
\label{cor:rigidity}
If the group $HY^2(\psi, \psi)$ is trivial, then $\psi$ is rigid.
\end{cor}


\section{Extending $2$-cocycles to deformations}
\label{sec:main}

The purpose of this section is to determine the obstructions for a $2$-cocycle in $CY^2(\psi,\psi)$ to be the infinitesimal of a deformation of $\psi$.  This is done by considering deformations modulo $t^N$ for $N = 1, 2, \ldots$ and determining the obstruction to extending a deformation modulo $t^N$ to a deformation modulo $t^{N+1}$.

\subsection{Deformations of finite order}
\label{subsec:def finite order}

Let $N$ be a positive integer.  A \emph{deformation of order $N$ of $\psi$} is simply a triple, $\Theta_t = (f_{D,t}; f_{E,t}; \Phi_t) = \sum_{i=0}^N t^i \theta_i$, which satisfies the conditions in Definition \ref{def:deformation psi} modulo $t^{N+1}$.  More explicitly, $f_{D,t} = \sum_{i=0}^N F_{D,i} t^i$ and $f_{E,t} = \sum_{i=0}^N F_{E,i} t^i$ satisfy equations $(5_\nu) - (9_\nu)$ in \cite[page 37]{mm} for $0 \leq \nu \leq N$, and $\Psi_t = \sum_{i=0}^N \psi_i t^i$ satisfies
   \begin{equation}
   \label{eq:psi f N original}
   \Psi_t\left(f^*_{D,t}(a,b)\right) ~\equiv~ f^*_{E,t}\left(\Psi_t(a), \Psi_t(b)\right) \quad \pmod{t^{N+1}},
   \end{equation}
or, equivalently,
   \begin{equation}
   \label{eq:psi f N}
   \sum_{i=0}^n \psi_i F^*_{D, n-i}(a, b) ~=~ \sum_{\substack{i + j + k = n \\ i, j, k \geq 0}} F^*_{E,i}\left(\psi_j(a), \psi_k(b)\right)
   \end{equation}
for $a, b \in D$, $0 \leq n \leq N$, and $* \in \lbrace l, r \rbrace$.  In particular, a deformation as defined in Section \ref{sec:deformation} can be regarded as a deformation of order $\infty$.

Given a deformation $\Theta_t$ of order $N$, it is said to \emph{extend to order $N + 1$} if and only if there exists a $2$-cochain $\theta_{N+1} = (F_{D,N+1}; F_{E,N+1}; \psi_{N+1}) \in CY^2(\psi,\psi)$ such that $\Thetabar_t = \Theta_t + t^{N+1} \theta_{N+1}$ is a deformation of order $N + 1$.  Such a $\Thetabar_t$ is called an \emph{order $N + 1$ extension of $\Theta_t$}.

Let $\Theta_t$ be a deformation of order $N$.  Consider the $3$-cochain 
   \begin{equation}
   \label{eq:def ob}
   \Ob(\Theta_t) ~=~ (\Ob_D; \, \Ob_E; \, \Ob_\psi) \in CY^3(\psi, \psi)
   \end{equation}
whose first two components are $\Ob_D = \sum_{\substack{i + j = N + 1 \\ i, j > 0}} F_{D,i} \circ F_{D,j}$, $\Ob_E = \sum_{\substack{i + j = N + 1 \\ i, j > 0}} F_{E,i} \circ F_{E,j}$, where $\circ$ is the pre-Lie product in $CY^*(D,D)$ or $CY^*(E,E)$ \cite[Definition 6.7]{mm}.  The last component $\Ob_\psi \in CY^2(D,E)$ is given by (for $a, b \in D$)
   \begin{subequations}
   \label{eq:ob psi}
   \begin{align}
   \Ob_\psi(\lbrack 2 1\rbrack \otimes (a, b)) & ~=~ \sumprime F^l_{E,i}(\psi_j(a), \psi_k(b)) ~-~ \sum_{i=1}^N \psi_i F^l_{D,N+1-i}(a, b), \label{eq:ob psi 1} \\
   \Ob_\psi(\lbrack 1 2\rbrack \otimes (a, b)) & ~=~ \sumprime F^r_{E,i}(\psi_j(a), \psi_k(b)) ~-~ \sum_{i=1}^N \psi_i F^r_{D,N+1-i}(a, b), \label{eq:ob psi 2}
   \end{align}
   \end{subequations}
where
   \begin{equation}
   \label{eq:sum prime}
   \sumprime ~=~ \sum_{\substack{i + j = N + 1 \\ i, j > 0 \\ k = 0}} ~+~ \sum_{\substack{i + k = N + 1 \\ i, k > 0 \\ j = 0}} ~+~ \sum_{\substack{j + k = N + 1 \\ j, k > 0 \\ i = 0}} ~+~ \sum_{\substack{i + j + k = N + 1 \\ i, j, k > 0}}.
   \end{equation}
The $3$-cochain $\Ob(\Theta_t)$ is called the \emph{obstruction class} of $\Theta_t$.

\begin{lemma}
\label{lem:main}
The obstruction class $\Ob(\Theta_t)$ is a $3$-cocycle.
\end{lemma}

Since the proof of this Lemma is rather long, it is postponed until the end of this section.  Assuming this Lemma for the moment, here is the main result of this section.

\begin{thm}
\label{thm:main}
Let $\Theta_t$ be a deformation of order $N$ of $\psi$.  Then $\Theta_t$ extends to a deformation of order $N + 1$ if and only if the cohomology class of $\Ob(\Theta_t)$ vanishes.  More precisely, if $\theta_{N+1} = (F_{D,N+1}; F_{E,N+1}; \psi_{N+1}) \in CY^2(\psi, \psi)$ is a $2$-cochain, then $\Thetabar_t = \Theta_t + t^{N+1} \theta_{N+1}$ is an order $N + 1$ extension of $\Theta_t$ if and only if $\Ob(\Theta_t) = \delta \theta_{N+1}$.
\end{thm}

\begin{cor}
\label{cor:main}
If the group $HY^3(\psi, \psi)$ is trivial, then every $2$-cocycle in $CY^2(\psi, \psi)$ is the infinitesimal of some deformation.
\end{cor}

\begin{proof}[Proof of Theorem \ref{thm:main}]
In fact, $\Thetabar_t$ is a deformation of order $N + 1$ if and only if the following three statements hold:
   \begin{enumerate}
   \item $\fbar_{D,t} = f_{D,t} + t^{N+1}F_{D,N+1}$ satisfies $(5_\nu) - (9_\nu)$ in \cite{mm} for $0 \leq \nu \leq N + 1$.
   \item $\fbar_{E,t} = f_{E,t} + t^{N+1}F_{E,N+1}$ satisfies $(5_\nu) - (9_\nu)$ in \cite{mm} for $0 \leq \nu \leq N + 1$.
   \item \eqref{eq:psi f N} holds for $n = N + 1$.
   \end{enumerate}
It is shown in \cite{mm} that (1) is equivalent to $\Ob_D = \delta F_{D,N+1}$.  Similarly, (2) is equivalent to $\Ob_E = \delta F_{E,N+1}$.  On the other hand, \eqref{eq:psi f N} with $n = N + 1$ is equivalent to $(\psi F_{D, N+1} - F_{E,N+1} \psi - \delta \psi_{N+1}) = \Ob_\psi$.  In other words, $(1) - (3)$ are equivalent to $\Ob(\Theta_t) = \delta \theta_{N+1}$, as claimed.
\end{proof}

\begin{proof}[Proof of Lemma \ref{lem:main}]
This is a rather long argument with a lot of bookkeeping.  It is known that $\Ob_D$ and $\Ob_E$ are $3$-cocycles in $CY^*(D,D)$ and $CY^*(E,E)$, respectively \cite[Theorem 3.5]{mm}.  Thus, it remains to show that
   \begin{equation}
   \label{eq:main condition}
   \psi\Ob_D ~-~ \Ob_E\psi ~-~ \delta \Ob_\psi ~=~ 0
   \end{equation}
in $CY^3(D,E) = \Hom(k\lbrack Y_3 \rbrack \otimes D^{\otimes 3}, E)$.  We will concentrate on the $3$-tree $y = \lbrack 3 2 1 \rbrack$.  The other four cases are proved similarly.  So let $a, b, c \in D$.  We have
   \begin{equation}
   \label{eq1:delta ob psi}
   \begin{split}
   (\delta \Ob_\psi)&(\lbrack 3 2 1\rbrack \otimes (a, b, c)) \\
   &= \psi(a) \dashv \left\lbrack \sumprime F^l_{E,i}(\psi_j(b), \psi_k(c)) \,-\, \sum_{i=1}^N \psi_i F^l_{D,N+1-i}(b,c)\right\rbrack \\
   &\relphantom{} ~-~ \left\lbrack \sumprime F^l_{E,i}(\psi_j(a \dashv b), \psi_k(c)) \,-\, \sum_{i=1}^N \psi_i F^l_{D,N+1-i}(a \dashv b, c)\right\rbrack \\
   &\relphantom{} ~+~ \left\lbrack \sumprime F^l_{E,i}(\psi_j(a), \psi_k(b \dashv c)) \,-\, \sum_{i=1}^N \psi_i F^l_{D,N+1-i}(a, b \dashv c) \right\rbrack \\
   &\relphantom{} ~-~ \left\lbrack \sumprime F^l_{E,i}(\psi_j(a), \psi_k(b)) \,-\, \sum_{i=1}^N \psi_i F^l_{D,N+1-i}(a, b) \right\rbrack \dashv \psi(c).
   \end{split}
   \end{equation}
In order to show that this is equal to $(\psi\Ob_D - \Ob_E\psi)(\lbrack 3 2 1\rbrack \otimes (a, b, c))$, we need to analyze every sum in it.

We begin with the third sum in \eqref{eq1:delta ob psi}.  It follows from \eqref{eq:psi f N} that, for each $j$, we have
   \begin{equation}
   \label{eq:psi j a b}
   \psi_j(a \dashv b) ~=~ \sum_{\substack{\alpha + \beta + \gamma = j \\ \alpha, \beta, \gamma \geq 0}} F^l_{E,\alpha}(\psi_\beta(a), \psi_\gamma(b)) ~-~ \sum_{\substack{\lambda + \mu = j \\ 1\leq \mu \leq j}} \psi_\lambda F^l_{D,\mu}(a,b). 
   \end{equation}
Substituting this into the third sum in \eqref{eq1:delta ob psi}, we can rewrite it as
   \begin{multline}
   \label{eq:3 delta ob}
   -\sumprime F^l_{E,i}(\psi_j(a \dashv b), \psi_k(c))
   ~=~ -\sumprime_{\substack{\alpha + \beta + \gamma = j \\ \alpha, \beta, \gamma \geq 0}} F^l_{E,i}(F^l_{E,\alpha}(\psi_\beta(a), \psi_\gamma(b)), \psi_k(c)) \\
    +~ \sumprime_{\substack{\lambda + \mu = j \\ 1 \leq \mu \leq j}} F^l_{E,i}(\psi_\lambda F^l_{D,\mu}(a,b), \psi_k(c)).
   \end{multline}
Here the first sum on the right-hand side (henceforth abbreviated as r.h.s.) is given by 
   \begin{equation}
   \label{eq:3a}
   \sumprime_{\substack{\alpha + \beta + \gamma = j \\ \alpha, \beta, \gamma \geq 0}} 
   ~=~ \sum_{\substack{i + \alpha + \beta + \gamma = N + 1 \\ i,\, \alpha + \beta + \gamma > 0 \\ k = 0 \\ \alpha, \beta, \gamma \geq 0}} 
   ~+~ \sum_{\substack{i + k = N + 1 \\ i, k > 0 \\ \alpha = \beta = \gamma = 0}} 
   ~+~ \sum_{\substack{\alpha + \beta + \gamma + k = N + 1 \\ k,\, \alpha + \beta + \gamma > 0 \\ i = 0 \\ \alpha, \beta, \gamma \geq 0}} 
   ~+~ \sum_{\substack{i + \alpha + \beta + \gamma + k = N + 1 \\ i, k,\, \alpha + \beta + \gamma > 0 \\ \alpha, \beta, \gamma \geq 0}}.
   \end{equation}
The second sum $\sumprime_{\substack{\lambda + \mu = j \\ 1 \leq \mu \leq j}}$ is obtained from $\sumprime$ \eqref{eq:sum prime} in a similar way by imposing the additional conditions, $\lambda + \mu = j$, $1 \leq \mu \leq j$.  More precisely, we have  
   \begin{equation}
   \label{eq:3b}
   \sumprime_{\substack{\lambda + \mu = j \\ 1 \leq \mu \leq j}}
   ~=~ \sum_{\substack{i + \lambda + \mu = N + 1 \\ i, \mu > 0;\, \lambda \geq 0 \\ k = 0}} 
   ~+~ \sum_{\substack{\lambda + \mu + k = N + 1 \\ \mu,\, k > 0;\, \lambda \geq 0 \\ i = 0}} 
   ~+~ \sum_{\substack{i + \lambda + \mu + k = N + 1 \\ i, k, \mu > 0;\, \lambda \geq 0 }}.
   \end{equation}
The same remarks apply below when we encounter such a construction again.  In particular, the first sum on the r.h.s.\ of \eqref{eq:3 delta ob} is the sum of four terms, corresponding to the four sums on the r.h.s.\ of \eqref{eq:3a}.  The first one of these four terms splits into a sum,
   \begin{equation}
   \label{eq:3a1}
   \begin{split}
   -\sum_{\substack{i + \alpha + \beta + \gamma = N + 1 \\ i,\, \alpha + \beta + \gamma > 0 \\ k = 0 \\ \alpha, \beta, \gamma \geq 0}} F^l_{E,i}(F^l_{E,\alpha}&(\psi_\beta(a), \psi_\gamma(b)), \psi_k(c)) \\
   &= -\sum_{\substack{i + \alpha = N + 1 \\ i, \alpha > 0}} F^l_{E,i}(F^l_{E,\alpha}(\psi(a), \psi(b)), \psi(c)) \\
   &\relphantom{} ~-~ \sum_{\substack{i + \alpha + \beta + \gamma = N + 1 \\ i,\, \beta + \gamma > 0 \\ \alpha, \beta, \gamma \geq 0}} F^l_{E,i}(F^l_{E,\alpha}(\psi_\beta(a), \psi_\gamma(b)), \psi(c)).
   \end{split}
   \end{equation}
Observe that the first term on the r.h.s.\ of \eqref{eq:3a1} is one of the two summands of $-(\Ob_E\psi)(\lbrack 3 2 1 \rbrack \otimes (a, b, c))$.

By applying a similar argument to the fifth term in \eqref{eq1:delta ob psi}, using \eqref{eq:psi f N} on $\psi_k(b \dashv c)$, one can rewrite it as 
   \begin{multline}
   \label{eq:5 delta ob}
   \sumprime F^l_{E,i}(\psi_j(a), \psi_k(b \dashv c)) ~=~ \sumprime_{\substack{\alpha + \beta + \gamma = k \\ \alpha, \beta, \gamma \geq 0}} F^l_{E,i}(\psi_j(a), F^l_{E,\alpha}(\psi_\beta(b), \psi_\gamma(c)) \\
   ~-~ \sumprime_{\substack{\lambda + \mu = k \\ 1 \leq \mu \leq k}} F^l_{E,i}(\psi_j(a), \psi_\lambda F^l_{D,\mu}(b, c)).
   \end{multline}
Just as above, the first term on the r.h.s.\ of \eqref{eq:5 delta ob} is a sum of four terms, similar to \eqref{eq:3a} except that the roles of $j$ and $k$ are interchanged.  One of these four terms is
   \begin{equation}
   \label{eq:5a2}
   \begin{split}
   \sum_{\substack{i + \alpha + \beta + \gamma = N + 1 \\ i,\, \alpha + \beta + \gamma > 0 \\ j = 0 \\ \alpha, \beta, \gamma \geq 0}} F^l_{E,i}(\psi_j(a)&, F^l_{E,\alpha}(\psi_\beta(b), \psi_\gamma(c))) \\
   &= \sum_{\substack{i + \alpha = N + 1 \\ i, \alpha > 0}} F^l_{E,i}(\psi(a), F^l_{E,\alpha}(\psi(b), \psi(c))) \\
   &\relphantom{} ~+~ \sum_{\substack{i + \alpha + \beta + \gamma = N + 1 \\ i,\, \beta + \gamma > 0 \\ \alpha, \beta, \gamma \geq 0}} F^l_{E,i}(\psi(a), F^l_{E,\alpha}(\psi_\beta(b), \psi_\gamma(c))).
   \end{split}
   \end{equation}
Observe that the first term on the r.h.s.\ of \eqref{eq:5a2} is the other summand of $-(\Ob_E\psi)(\lbrack 3 2 1 \rbrack \otimes (a, b, c))$.

Now we consider the fourth term in \eqref{eq1:delta ob psi}.  For each $i = 1, \ldots , N$, we use equation $(5_{N+1-i})$ in \cite{mm} to obtain
   \begin{multline}
   \label{eq1:4}
   F^l_{D,N+1-i}(a \dashv b, c) ~=~ \sum_{j=0}^{N+1-i} F^l_{D,j}(a, F^l_{D,N+1-i-j}(b,c)) \\ ~-~ \sum_{j=0}^{N-i} F^l_{D,j}(F^l_{D,N+1-i-j}(a,b),c).
   \end{multline}
Substituting this into the fourth term on the r.h.s.\ of \eqref{eq1:delta ob psi}, it can be rewritten as
   \begin{equation}
   \label{eq2:4}
   \begin{split}
   \sum_{i=1}^N \psi_i F^l_{D,N+1-i}(a \dashv b, c) &~=~ \sum_{i=1}^N \psi_i F^l_{D,N+1-i}(a, b \dashv c) \\
   &\relphantom{} ~+~ \sum_{\substack{i + j + k = N + 1 \\ i, k > 0;\, j \geq 0}} \psi_i F^l_{D,j}(a, F^l_{D,k}(b,c)) \\
   &\relphantom{} ~-~ \sum_{\substack{i + j + k = N + 1 \\ i, k > 0;\, j \geq 0}} \psi_i F^l_{D,j}(F^l_{D,k} (a, b), c).
   \end{split}
   \end{equation}
Observe that the first term on the r.h.s.\ of \eqref{eq2:4} cancels with the sixth term on the r.h.s.\ of \eqref{eq1:delta ob psi}.

On the other hand, the second term on the r.h.s.\ of \eqref{eq2:4} can be expanded as
   \begin{multline}
   \label{eq1:4b} 
   \sum_{\substack{i + j + k = N + 1 \\ i, k > 0;\, j \geq 0}} \psi_i F^l_{D,j}(a, F^l_{D,k}(b,c)) 
   ~=~ \sum_{k=1}^N \left\lbrack \sum_{\substack{i + j = N + 1 - k \\ i, j \geq 0}} \psi_i F^l_{D,j}(a, F^l_{D,k}(b, c))\right\rbrack \\
   ~-~ \psi \left\lbrack \sum_{\substack{j + k = N + 1 \\ j, k > 0}} F^l_{D,j}(a, F^l_{D,k}(b,c)) \right\rbrack.
   \end{multline}
Observe that the second term on the r.h.s.\ of \eqref{eq1:4b} is one of the two summands of $(\psi \Ob_D)(\lbrack 3 2 1 \rbrack \otimes (a, b, c))$.  For $1 \leq k \leq N$, one has $1 \leq N + 1 - k \leq N$, and so \eqref{eq:psi f N} allows us to rewrite the first term on the r.h.s.\ of \eqref{eq1:4b} as
   \begin{equation}
   \label{eq:4b1}
   \begin{split}
   \sum_{k=1}^N & \left\lbrack\sum_{\substack{i + j = N + 1 - k \\ i, j \geq 0}} \psi_i F^l_{D,j}(a, F^l_{D,k}(b, c))\right\rbrack \\
   &=~ \sum_{k=1}^N \left\lbrack \sum_{\substack{\alpha + \beta + \gamma = N + 1 - k \\ \alpha, \beta, \gamma \geq 0}} F^l_{E,\alpha}(\psi_\beta(a), \psi_\gamma F^l_{D,k}(b, c))\right\rbrack \\
   &=~ \psi(a) \dashv \left\lbrack \sum_{i=1}^N \psi_i F^l_{D,N+1-i}(b,c) \right\rbrack ~+~ \sumprime_{\substack{\lambda + \mu = k \\ 1 \leq \mu \leq k}} F^l_{E,i}(\psi_j(a), \psi_\lambda F^l_{D,\mu}(b, c)).
   \end{split}
   \end{equation}  
In the last line, the two terms cancel with the second terms on the r.h.s.\ of, respectively, \eqref{eq1:delta ob psi} and \eqref{eq:5 delta ob}.

A similar argument, applied to the last term in \eqref{eq2:4}, yields
   \begin{equation}
   \label{eq:4c}
   \begin{split}
   -&\sum_{\substack{i + j + k = N + 1 \\ i, k > 0;\, j \geq 0}} \psi_i F^l_{D,j}(F^l_{D,k} (a, b), c) \\
   &~=~ -\left\lbrack \sum_{i=1}^N \psi_i F^l_{D,N+1-i}(a, b) \right\rbrack \dashv \psi(c) 
   ~-~ \sumprime_{\substack{\lambda + \mu = j \\ 1 \leq \mu \leq j}} F^l_{E,i}(\psi_\lambda F^l_{D,\mu}(a,b), \psi_k(c)) \\
   &\relphantom{} ~+~ \psi \left\lbrack \sum_{\substack{j + k = N + 1 \\ j, k > 0}} F^l_{D,j}(F^l_{D,k}(a, b),c) \right\rbrack.
   \end{split}
   \end{equation}
On the r.h.s.\ of \eqref{eq:4c}, the last term is the other summand of $(\psi \Ob_D)(\lbrack 3 2 1 \rbrack \otimes (a, b, c))$, while the first two terms cancel with the last terms on the r.h.s.\ of, respectively, \eqref{eq1:delta ob psi} and \eqref{eq:3 delta ob}.

The argument so far tells us that the following element in $E$,
   \begin{equation}
   \label{eq2:delta ob psi}
   (\delta \Ob_\psi - \psi\Ob_D + \Ob_E\psi)(\lbrack 3 2 1 \rbrack \otimes (a, b, c)),
   \end{equation}
is equal to the following sum (where $\alpha, \beta, \gamma \geq 0$ wherever applicable):
   \begin{equation}
   \label{eq:monster sum}
   \begin{split}
   &\sum_{\substack{i+j = N+1 \\ i,j>0}} F^l_{E,i}(\psi_j(a), \psi(b) \dashv \psi(c)) 
   ~+~ \sum_{\substack{i + \alpha + \beta + \gamma = N + 1 \\ i,\, \beta + \gamma > 0}} F^l_{E,i}(\psi(a), F^l_{E,\alpha}(\psi_\beta(b), \psi_\gamma(c))) \\
   &~+~ \sum_{\substack{j + \alpha + \beta + \gamma = N+1 \\ 1 \leq j \leq N}} \psi_j(a) \dashv F^l_{E,\alpha}(\psi_\beta(b), \psi_\gamma(c)) \\
   &~+~ \sum_{\substack{i + j + \alpha + \beta + \gamma = N+1 \\ i, j, \, \alpha + \beta + \gamma > 0 }} F^l_{E,i}(\psi_j(a), F^l_{E,\alpha}(\psi_\beta(b), \psi_\gamma(c))) \\
   &~-~ \sum_{\substack{i + \alpha + \beta + \gamma = N + 1 \\ i,\, \beta + \gamma > 0}} F^l_{E,i}(F^l_{E,\alpha}(\psi_\beta(a), \psi_\gamma(b)), \psi(c)) \\
   &~-~ \sum_{\substack{i+k = N+1 \\ i,k>0}} F^l_{E,i}(\psi(a) \dashv \psi(b), \psi_k(c)) \\
   &~-~ \sum_{\substack{\alpha + \beta + \gamma + k = N + 1 \\ k,\, \alpha + \beta + \gamma > 0}} F^l_{E,\alpha}(\psi_\beta(a), \psi_\gamma(b)) \dashv \psi_k(c) \\
   &~-~ \sum_{\substack{i + \alpha + \beta + \gamma + k = N + 1 \\ i, k, \, \alpha + \beta + \gamma > 0}} F^l_{E,i}(F^l_{E,\alpha}(\psi_\beta(a), \psi_\gamma(b)), \psi_k(c)) \\
   &~+~ \left\lbrack\psi(a) \dashv \sumprime F^l_{E,i}(\psi_j(b), \psi_k(c))\right\rbrack 
   ~-~ \left\lbrack \sumprime F^l_{E,i}(\psi_j(a), \psi_k(b)) \right\rbrack \dashv \psi(c).
   \end{split}
   \end{equation} 
This sum can be written more compactly as
   \begin{equation}
   \label{eq:monster sum'}
   \widetilde{\sum} \left\lbrack F^l_{E,\lambda}(\psi_\alpha(a), F^l_{E,\mu}(\psi_\beta(b), \psi_\gamma(c))) ~-~ F^l_{E,\lambda}(F^l_{E,\mu}(\psi_\alpha(a), \psi_\beta(b)), \psi_\gamma(c)) \right\rbrack,
   \end{equation}
where
   \begin{equation}
   \label{eq:sumtilde}
   \widetilde{\sum} 
   ~=~ \sum_{\substack{\alpha + \beta + \gamma + \lambda + \mu = N + 1 \\ 1 \leq \alpha + \beta + \gamma \leq N \\ \alpha, \beta, \gamma, \lambda, \mu \geq 0}} 
   ~+~ \sum_{\substack{\alpha + \beta = N + 1 \\ \alpha, \beta > 0 \\ \lambda, \mu, \gamma = 0}} 
   ~+~ \sum_{\substack{\alpha + \gamma = N + 1 \\ \alpha, \gamma > 0 \\ \lambda, \mu, \beta = 0}}
   ~+~ \sum_{\substack{\beta + \gamma = N + 1 \\ \beta, \gamma > 0 \\ \lambda, \mu, \alpha = 0}}
   ~+~ \sum_{\substack{\alpha + \beta + \gamma = N + 1 \\ \alpha, \beta, \gamma > 0 \\ \lambda, \mu = 0}}.
   \end{equation}
In particular, it follows from the expression \eqref{eq:monster sum'}, equation $(5_{\lambda + \mu})$ in \cite{mm}, and one of the dialgebra axioms (the associativity of $\dashv$) that the sum in \eqref{eq:monster sum}, and hence $(\delta \Ob_\psi - \psi\Ob_D + \Ob_E\psi)(\lbrack 3 2 1 \rbrack \otimes (a, b, c))$, is equal to $0$.

The argument that $(\delta \Ob_\psi - \psi\Ob_D + \Ob_E\psi)(y \otimes (a, b, c))$ is equal to $0$ for the other four $3$-trees $y \in Y_3$ is similar to the one given above.  Instead of equation $(5_\nu)$ in \cite{mm} and the associativity of $\dashv$, one makes use of $(6_\nu)$, $(7_\nu)$, $(8_\nu)$, or $(9_\nu)$ in conjunction with one of the other four dialgebra axioms.

This finishes the proof of Lemma \ref{lem:main}.
\end{proof}



\begin{thebibliography}{99}
\bibitem{fra}A.\ Frabetti,  Dialgebra (co)homology with coefficients, in:  Dialgebras and related operads, 67-103, Lecture Notes in Math., 1763, Springer, Berlin, 2001. 

\bibitem{freiger}Y. Fr\'eiger, A new cohomology theory associated to deformations of Lie algebra morphisms, Lett. Math. Phys. 70 (2004), 97-107.

\bibitem{ger0}M.\ Gerstenhaber, The cohomology structure of an associative ring, Ann.\ Math.\ \textbf{78} (1963), 267-288.

\bibitem{ger1}M.\ Gerstenhaber, On the deformation of rings and algebras, Ann.\ Math.\ \textbf{79} (1964), 59-103.

\bibitem{gs1}M.\ Gerstenhaber and S.\ D.\ Schack, On the deformation of algebra morphisms and diagrams, Trans. Amer. Math. Soc. 279 (1983), 1-50. 

\bibitem{gs2}M.\ Gerstenhaber and S.\ D.\ Schack, On the cohomology of an algebra morphism,  J. Algebra 95 (1985), 245-262.

\bibitem{gs3}M.\ Gerstenhaber and S.\ D.\ Schack, Sometimes $H\sp 1$ is $H\sp 2$ and discrete groups deform, in: Geometry of group representations (Boulder, CO, 1987),  149-168, Contemp. Math., 74, Amer. Math. Soc., Providence, RI, 1988.

\bibitem{loday}J.-L.\ Loday, Dialgebras, in:  Dialgebras and related operads, 7-66, Lecture Notes in Math., 1763, Springer, Berlin, 2001.

\bibitem{lp}J.-L.\ Loday and T.\ Pirashvili,  Universal enveloping algebras of Leibniz algebras and (co)homology, 
Math.\ Ann.\ 296 (1993), 139-158.

\bibitem{mm}A.\ Majumdar and G.\ Mukherjee, Deformation theory of dialgebras, K-theory 27 (2002), 33-60.

\bibitem{nr}A.\ Nijenhuis and R.\ W.\ Richardson, Deformations of homomorphisms of Lie algebras, Bull. Amer. Math. Soc. 73 (1967), 175-179.

\end{thebibliography}
\end{document}